\font\titlefont =cmr17 
\font \didesnisfont cmbx12 
\font\eightit= cmti8 
%
\hsize=12cm
\vsize=20 cm
\maxdepth=2.2pt
\hoffset=1cm
\voffset=1cm
\newdimen\tarpas  \tarpas= 10pt
\newdimen\isvirsaus \isvirsaus=5pt

\newcount\propositionnumber 
\propositionnumber=0
\newtoks\maintoks
\newtoks\indextoks
%
\newcount\skaitliukas 

\def\literalnumber #1*{
\ifnum #1<10 \hskip 1em  #1
\else \ifnum #1< 100 \hskip 0.5em  #1
\else #1
\fi\fi}
%
\long\def\Skyr #1.#2.. {}%

\long\def\Proposition  #1\Proof:{%
\advance \propositionnumber by 1
\medbreak\noindent {\bf Proposition
\the\propositionnumber.\enspace}
{\sl #1}\par
\ifdim\lastskip<\medskipamount\removelastskip\penalty55\medskip\fi
\noindent{\bf Proof: }
}
%
%

\def\sqr#1#2{{\vcenter{\vbox{\hrule height.#2pt
\hbox{\vrule width.#2pt height#1pt \kern#1pt
\vrule width.#2pt}\hrule height.#2pt}}}}



\def\baf #1*{{\bf #1}}
%
\def\virsutskaitlius{%
\vbox{
\ifodd\pageno 
\rightheadline
\else 
\leftheadline
\fi
\vss}
}

\def\leftheadline{}
\def\rightheadline{}
\def\apatskaitlius{%
\vbox{
\ifnum \pageno<2 
\else 
\ifodd\pageno 
\rightnumber
\else 
\leftnumber
\fi
\fi
\vss}
}

\newdimen\ruleht
\ruleht=.3pt 
\def\leftnumber{\hbox to \hsize{\vbox to 8pt{}
{\eightit\hskip 1pt \the\pageno\hfil }}}
\def\rightnumber{\hbox to \hsize{\vbox to 8pt{}
{\eightit\hfil\the\pageno }}}
\pageno=1

\def\viengubasspausdinimas{%
\shipout\vbox{%
\vskip\isvirsaus
\vbox{\virsutskaitlius\vskip\isvirsaus\box255\apatskaitlius}
\vskip\isvirsaus
}
\global\advance \pageno by 1
}

\def\Rym#1'{\uppercase\expandafter{\romannumeral#1}}%
%
\def\specialtilde{\ifmmode
\catcode `\^=7 
\def\next{^}\else\let\next\normalref%
\fi\next}

\catcode`\^=\active\let ^=\specialtilde %
\catcode`\~=\active \def~{\penalty100000} 

\newtoks\maint
%

\def\normalref#1{%
\atskiria#1\gl%
{\it
\pakeistas}
}
\def\atskiria#1=#2\gl{
\edef\pagrindinis{\if=#2\else 
#2\fi}
\def\pakeistas{#1}%
\if=#2\else \def\pakeistas{#1#2}\fi%
\indextoks{{#1\it \pagrindinis}}
}
\def\vect#1 {\setbox0=\hbox{\.#1}%
\vbox{\hbox to \wd0{\hss $\rightarrow$\hss}%
\vskip -4pt \copy0} }%
\def\leaderfill{\leaders\hbox to 1em{\hss.\hss}\hfill}
\def\fract#1:#2 {{#1\over#2}}




\def\bare#1 {\setbox0=\hbox{#1}%
\vbox{\hbox to \wd0{\hss $-$\hss}%
\vskip -4pt \copy0} }%

\mathchardef\tustis="001F 

%
%

%
\def\maprightdviguba #1 #2 {
\raise 1pt\hbox{$\mapright {\lower 1pt\hbox{$\scriptstyle #1$}} $}%
\llap{\lower 1pt\hbox{$\maprightap {\raise 1pt \hbox{$\scriptstyle #2$}} $}}}

\def\mapleftright #1 #2 {
\raise 1pt\hbox{$\mapleft {\lower 1pt\hbox{$\scriptstyle #1$}} $}%
\llap{\lower 1pt\hbox{$\maprightap {\raise 1pt \hbox{$\scriptstyle #2$}} $}}}

\def\mapdownk #1 {\mathrel{#1\kern 2pt \downarrow}}
\def\mapdownd #1 {\mathrel{\downarrow \kern 2pt #1}}
\def\mapupk #1 {\mathrel{#1\kern 2pt \uparrow}}
\def\mapupd #1 {\mathrel{\uparrow \kern 2pt #1}}
\def\mapright #1 {\mathrel{\smash{\mathop{\longrightarrow}\limits^{#1}}%
\vphantom\longrightarrow}}
\def\maprightap #1 {\mathrel{\smash{\mathop{\longrightarrow}\limits_{#1}}%
\vphantom\longrightarrow}}
\def\mapleft #1 {\mathrel{\smash{\mathop{\longleftarrow}\limits^{#1}}%
\vphantom\longleftarrow}}
\def\mapleftap #1 {\mathrel{\smash{\mathop{\longleftarrow}\limits_{#1}}%
\vphantom\longleftarrow}}
%
\def\tworightarrows{\mathop{\raise 1pt \hbox{$\rightarrow$}\llap{\raise -1pt
\hbox{$\rightarrow$}}}}
\def\rightleftarrows{\mathop{\raise 1pt \hbox{$\rightarrow$}\llap{\raise -1pt
\hbox{$\leftarrow$}}}}
%
%
%
%
%
%
%
\def\lim{\mathop{\tt lim}}%

\def\CAT{\mathord{{\bf CAT}}}%
\def\SET{\mathord{\bf SET}}%
\def\GRPH{\mathord{\bf GRPH}}%

\def\and{\mathord {{\tt and }}}
'

%

%
%
%
%
%
%
%
%
%


%
%
%
%
%

%
%
%
%
%
%
%
%
%
%
%
%
%
%
%
%
%
%
%
%

%
%

%

%
%

\mathchardef\tustis="001F 
\mathchardef\bangele="0365 

\def\onetilde#1 {\setbox 0=\hbox{$#1$}
\mathord {\vbox {\hbox to \wd0{$\bangele$}\vskip -9.5pt \copy 0}}}

\def\twotilde #1 {\setbox 0=\hbox{$#1$}
\mathord{ \vbox{\vskip -2pt \hbox to \wd0{ \hss $\bangele$ \hss}
\vskip -9.5 pt  \hbox to \wd0{ \hss $\bangele$ \hss}\vskip -9.5 pt
  \copy 0}}}
\def\threetilde #1 {\setbox 0=\hbox{$#1$}
\mathord{ \vbox{\vskip -5pt \hbox to \wd0{ \hss $\bangele$ \hss}\vskip -9.5 pt
 \hbox to \wd0{ \hss $\bangele$ \hss}\vskip -9.5 pt
 \hbox to \wd0{ \hss $\bangele$ \hss}\vskip -9.5 pt
  \copy 0}}}
%
%
%



%
%

\vskip 1cm
\centerline{\titlefont  Closed categories}
{\didesnisfont Gintaras Valiukevi\v cius \hfil  29 June 2022 }

{\it Veja

 Ukmerg\.e, Lithuania

e-mail: arklys@vivaldi.net
}
\vskip\tarpas
{\bf Abstract.}
We are checking the closed categories beginning with the category of all
sets $\SET$ and ending with the category of all categories $\CAT$.
The novelty is a generalizing the adjoint functors to the category of
directed graphs $\GRPH$. We have described  for what condition we get
bijective name mapping for graphs transports
$\GRPH(X\times T;Y)\rightarrow \GRPH (X;Y^T)$.
As an example of future applications we introduce a notion of relator
to some tensor product.

\vskip\tarpas
\def \.{./flash/pmetric}
\input./closed.tex
\vskip\tarpas
 REFERENCES
\vskip \tarpas

1. Stefan Banach 1922, Sur les op\'eration  dans les ensembles abstraits
et leurs applications  aux \'equations int\'egrales, 7-33,
Fund. Math. {\bf 3}.

2. J. Schauder 1927. Zur Theorie stetiger Abbildungen  in Funktionalr\"aumen,
47 -65, 417 - 431, Math. Z. {\bf 26}.

3. S. Karlin 1948, Bases in Banach spaces, 971 - 985, Duke Math. J. {\bf 15}.

4, Nicolas Bourbaki, Th\'eorie des ensembles ch. 4, Hermann :Paris 1957,
Mir: Moskva 1965.

5. Nicolas Bourbaki, Alg\`ebre chap 1-3, Masson: Paris 1970,  Chap.4-6,\break
Hermann: Paris, Nauka: Moskva 1965.

6. Jean Benabou 1963, Cat\'egories avec multiplication, 1887 - 1890, C. R. Paris
{\bf 256}.

7. A. Pietsch 1963, Zur Theorie der topologischen Tensorprodukte, 19 - 30, Math.
Nachr. {\bf 25}.

8. Helmut Schaefer 1966, Topological vector spaces, MacMillan: N. Y. , \break
London,
1971 Mir: Moskva.

9 Saunders MacLane 1971, Categories for the working mathematician,\break
Sprin\-ger: New York/ Berlin/ Heidelberg,  1997 second edition.

10. G. Valiukevi\v cius 1992, Integration in ordered spaces, I, Veja: Vilnius.

11. R. Backhouse, P. Hoogedijk 2000, Elements of a relational theory of data\-types.

12. G. Valiukevi\v cius 2009, United sight to an algebraic operations and
convergence, arXiv.org.

\vskip\tarpas
21 Septembre 2022.
\end